\documentclass{article}

\usepackage{amsmath}
\usepackage{amssymb}
\usepackage{amsfonts}
\usepackage{amsopn}
\usepackage[all]{xy}
\usepackage{amsthm}
\usepackage{amscd}

\swapnumbers
\theoremstyle{plain}
\newtheorem{thm}{Theorem}[section]
\newtheorem{lemma}[thm]{Lemma}
\newtheorem{sublemma}[thm]{Sublemma}

\theoremstyle{definition}
\newtheorem{rem}[thm]{Remark}
\newtheorem{ntt}[thm]{ }

\newcommand{\ra}{\rightarrow}

\newcommand{\AS}{{\mathbb A}}

\newcommand{\Bb}{{\mathcal B}}

\title{On Knebusch's Norm Principle for quadratic forms over semi-local rings}
\author{K.~Zainoulline}

\begin{document}

\maketitle

\begin{abstract}
We prove the version of Knebusch's
Norm Principle for simple extensions of (semi-)local rings.
As an application we prove
the Grothendieck-Serre's conjecture 
on principal homogeneous spaces for the split case of the spinor group.
\end{abstract}

\section{Introduction}

Let $L/F$ be a finite field extension and $q$ be a regular quadratic
form over $F$. The well-known Norm Principle for quadratic forms 
due to M.~Knebusch
says that $N(qL)\subset qF$, where $N$ is the norm map and $qL$ is
the group of values of quadratic form $q$ over $L$.

In the present paper we prove the version of Knebusch's
Norm Principle for simple extensions of (semi-)local rings
(the extension $S/R$ is called simple, where $R$ is local, 
if $S=R[t]/(p(t))$ is 
the quotient of the polynomial ring $R[t]$ by some monic polynomial $p$).
A finite etale extension (or separable in the case of fields)
gives an example of a simple extension.
This result is the key step in the proof of 
Grothendieck-Serre's conjecture
on principal homogeneous space for the case of of the spinor group 
(Theorem \ref{groth}):
the conjecture says that
the canonical map $H_{et}^1(R,Spin_q)\ra H^1_{et}(K,Spin_q)$
has the trivial kernel, where $R$ is a local regular ring and
$K$ is it's quotient field.

The main result is Theorem \ref{normpr} which says the Norm Principle
holds if $R$ is a local domain with the residue field of characteristic 0.
Briefly speaking, we use the fact that any square is represented
by the product of the values of quadratic  form $q$. 
Hence, we may consider the original arguments (see \cite{La}) ``modulo
squares''.
This together with some  general position arguments (see \ref{genpos}) 
allows us to control the leading coefficient of the polynomial $h$
(see the equation (**)) in order to make it invertible
(which is the key point in applying the induction hypotheses).

\ 

The author is grateful to Max-Planck-Institute f\"ur Mathematik and European
Post-Doctoral Institute
for hospitality and financial support.

\section{Simple extensions of local rings}
\begin{ntt}
All rings are assumed to be commutative with units.
By $R$ we denote a local domain with the residue field $k$.
By ``bar'' we mean the reduction modulo the maximal ideal, i.e., 
$k=\bar{R}$.
By a ring extension $S/R$ we mean a ring $S$ together
with a ring monomorphism $R\ra S$.
By $S^*$ we denote the group of invertible elements of $S$.

Let $S/R$ be a ring extension such that $S$ is free as the $R$-module.
Then there is the norm map denoted by $N^S_R$ that is given as follows.
For any $c\in S$ let $l_c:S\ra S$ be the endomorphism of the free
$R$-module $S$ given by the left
multiplication by $c$. Then we set $N^S_R(c)=\det(l_c)$, where $l_c$
is the respective matrix.
\end{ntt}

\begin{ntt}
We say a ring extension $S/R$ is a simple extension of degree $n$ 
if there exists an element $c\in S^*$ such that $S$ is free as the $R$-module
with the basis $\{1, c,\ldots , c^{n-1}\}$ denoted by $\Bb(c)$. The element $c$
is called a primitive element.
Let $c^n=a_0+a_1c+\ldots+ a_{n-1}c^{n-1}$ be the unique presentation
of the element $c^n$ in the basis $\Bb(c)$.
The polynomial $p_c(t)=t^n-a_{n-1}t^{n-1}-\ldots -a_0$ is called
a minimal polynomial for the primitive element $c$.
It has the property $p_c(c)=0$. Observe that $S$ can be identified
with the quotient $R[t]/p_c(t)$ of the polynomial ring $R[t]$
modulo the principal ideal generated by $p_c(t)$.
Since $c$ is invertible, it's norm $N^S_R(c)=(-1)^{n-1}p_c(0)=(-1)^na_0$ 
is invertible as well.

To the opposite direction, for a given monic 
polynomial $p(t)\in R[t]$ such that $p(0)\in R^*$ 
the ring extension $S=R[t]/p(t)$ over $R$ 
is a simple extension of degree $n=\deg p$.
The image of $t$ by means of the canonical map $R[t]\ra R[t]/p(t)$
gives the respective primitive element of the extension $S/R$. 
\end{ntt}

\begin{ntt}
Clearly, if $c\in S^*$ is a primitive element of 
a simple extension $S/R$,
then it's inverse $c^{-1}$ is primitive and
$r_1c+r_0$ is primitive for any $r_1\in R^*$ and 
$r_0\in R$ such that $p_c(-r_0/r_1)\in R^*$. 
\end{ntt}

\begin{ntt}\label{primelts}
Let $S/R$ be a simple extension with the primitive element $c$,
then $\bar{S}/\bar{R}$ is a simple extension of the same degree
with the primitive element $\bar{c}$. Moreover, if $b$ is a primitive
element of $\bar{S}$, then by Nakayama's lemma 
the preimage $c\in \rho^{-1}(b)$ is primitive as well,
where $\rho: S \ra \bar{S}$ is the reduction map. 
In other words, if $S^*_{prim}$ denotes the subset of primitive elements of 
$S$ and $\bar{S}^*_{prim}$ 
denotes the subset of primitive elements of $\bar{S}$,
then we have 
$S^*_{prim}=\rho^{-1}(\bar{S}^*_{prim})$. 
\end{ntt}

\begin{rem} Let $S/R$ be a simple extension.
In the case $R=k$ is a field the algebra $S$ can be
viewed as the product of local Artinian algebras over $k$.
For instance, if $k$ is algebraically closed, the algebra $S$ 
is isomorphic to the product of algebras of the kind $k[t]/t^m$, $m\geq 1$.

In the case $S$ is a field 
we get a simple field extension $S/k$. Recall 
(by the Primitive Element Theorem) that any finite separable field
extension is simple but not in the other direction. There are examples of
finite field extensions which are not simple.
\end{rem}

\section{Some general position arguments}

In the present section $S$ will be a simple extension
of an infinite field $k$.

\begin{ntt}\label{topol}
Let $S$ be a simple extension of $k$ of degree $n$.
The $k$-algebra $S$ can be viewed as the $k$-vector space of dimension $n$
and, thus, can be identified with the set of rational points of the
affine space $\AS^n$ over $k$. From this point on we assume
$S=\AS^n(k)$ is the topological space by means of Zariski topology structure.

For example, any map $S\ra S$ given by $b\mapsto f(b)$ is continuous,
where $f(t)\in S[t]$ is a polynomial with coefficients in $S$.
And the set $S^*$ of invertible elements 
is open in $S$, since it is given by the equation $N^S_k(x)\neq 0$,
where $N^S_k$ is the norm map.

In the case of an infinite field $k$ 
this topology has the important property 
-- the intersection
of any two open subsets is non-empty.
\end{ntt}

\begin{lemma}
Let $S$ be a simple extension of a field $k$. 
Then the subset of primitive elements
$S^*_{prim}$ is non-empty and open in $S^*$. 
\end{lemma}

\begin{proof}
We fix some primitive element $c$ of $S$. An element $b\in S^*$
is primitive iff the matrix $(b_{i,j})_{i,j=0}^{n-1}$ 
has non-zero determinant, where $b_{i,j}$ is the $i$-th
coefficient in the presentation of the $j$-th power $b^j$ of the element $b$
in the basis $\Bb(c)$. 
Observe that $b_{i,1}=b_i$ are the coefficients of the presentation 
of the element $b$ in the basis $\Bb(c)$, i.e,
$b=b_0+b_1c+\ldots + b_{n-1}c^{n-1}$.
Clearly, the determinant $\det(b_{i,j})$ can be viewed as the polynomial
in $n$ variables $b_0,\ldots,b_{n-1}$ and, thus,
the subset 
$$
S^*_{prim}=\{b=(b_0,\ldots,b_{n-1})\in S^*|\det(b_{i,j})\neq 0\}
$$ 
is open in $S^*$.
\end{proof}

\begin{lemma}\label{cprim}
Let $S$ be a simple extension of an infinite field $k$ 
of characteristic different from 2 and 
let $c\in S^*$ be an invertible element.
Then the subset 
$V_c=\{b\in S^*| cb^2\; \text{is primitive} \}$ is non-empty and open in $S^*$.
\end{lemma}

\begin{proof}
Since the map $S^*\ra S^*$ given by $b\mapsto cb^2$ is continuous, 
it is enough to show that the image of the map 
$f:S^*\ra S^*$ given by $b\mapsto b^2$ is dense.
(observe that the multiplication by $c$ is the homeomorphism).
The algebra $S$ splits as the product of local Artinian algebras and
the image of $f$ is dense if the image of the restriction
of $f$ to the each component of this product is dense.
Thus, we may assume $S$ is a local Artinian algebra over $k$ and,
hence, it is irreducible.
 
Assume the image of $f$ is not dense. 
Then the closure of the image of $f$ in $S$
must have the dimension strictly less than the dimension of $S$ (considered
as the affine space over $k$). It means that $f$ induces the regular
map between two affine spaces such that the dimension of the target
space is strictly less than the dimension of the origin space $S$.
In particular, there exists an element $u\in S^*$ such that
the equation $b^2=u$ has infinite number of solutions
(the dimension of the fiber of $f$ over $u$ is $\geq 1$ and $k$ is infinite).
Hence, we get contradiction by Lemma \ref{densit}. 
\end{proof}

\begin{lemma}\label{densit}
Let $S$ be a simple extension over an infinite field $k$
of characteristic different from 2. 
Let $u$ be an invertible element of $S$.
Then the number of solution of the equation $b^2=u$ is finite (or empty).
\end{lemma}

\begin{proof}
Let $k'$ be the algebraic closure of the field $k$.
Let $S'=S\otimes_k k'$ be the base change of $S$.
Clearly, $S'$ is the simple extension of $k'$ of the same degree.
The number of solutions of $b^2=u$ over $S$ is finite
if it is finite over $S'$. 

Since the algebra $S'$ splits as the finite product of algebras of the kind
$A_m=k'[t]/t^m$, $m\geq 1$, it is enough to show that the number
of solutions of $b^2=u$ is finite in $A_m$ for any $m$.

The case $m=1$ is trivial, since $A_1=k'$ is a field.
Let $m>1$.
In the basis $\Bb(t)$ of $A_m$ our equation can be written as:
$$
(b_0+b_1t+\ldots+b_{m-1}t^{m-1})^2=u_0+u_1t+\ldots + u_{m-1}t^{m-1}.
$$
Hence, we get the system of $m$ quadratic equations over $k'$
\begin{equation}
\tag{*}
b_0^2=u_0,\; 2b_0b_1=u_1,\; 2(b_0b_2+b_1^2)=u_2,\; 2(b_3b_0+b_2b_1)=u_3,\ldots
\end{equation}
which has the property that any element $b_j$ is the solution of
the quadratic or linear equation over $k'$ (precisely the $j+1$-th equation) 
with coefficients $b_i$, $i<j$, and $u_i$, $i\leq j$.
Then it follows immediately that the number of solutions of (*) is finite.
\end{proof}

\begin{rem}
The assumption that the characteristic of $k$ is different from 2 is essential.
Take $k$ to be the algebraic closure of the finite field ${\mathbb F}_2$
then the algebra $S=k[t]/t^2$ is the simple extension of $k$.
But it easy to see that the image of the map $b\ra b^2$ coincides
with the subspace $k\cdot 1$ in $S=k\cdot1 \oplus k\cdot t$ 
which consists of all non-primitive
elements of $S$.
\end{rem}

\begin{ntt}
Let $c\in S^*$ be a primitive element of 
a simple extension $S/k$ of degree $n$.
Let $x$ be an element of $S$. 
By the symbol $\{x,c\}$
we denote the $n$-th coefficient
of the presentation of $x$
in the basis $\Bb(c)$, i.e., $\{x,c\}=x_{n-1}$, where 
$\sum_{i=0}^{n-1} x_ic^i=x$.
\end{ntt}

\begin{lemma}\label{cxprim}
Let $S$ be a simple extension of degree $n$ of a field $k$ of
characteristic 0.
Let $c\in S^*_{prim}$ be a primitive element and let $x$ be a non-zero element
of $S$.
Then the subset $W_{c,x}=\{b\in V_c|\{xb^{-1},cb^2\}\neq 0\}$  
is non-empty and open in $S^*$, where $V_c$ is the non-empty open subset
from Lemma \ref{cprim}. 
\end{lemma}

\begin{rem}\label{charrem} The assumption that $k$ has characteristic 0 is essential.
Take $k$ to be the algebraic closure of the finite field ${\mathbb F}_3$.
Consider the algebra $S=k[t]/t^3-1$. Take $x=c=t$. It is easy to see
that $\{xb^{-1},cb^2\}=0$ for all $b\in S^*$.
\end{rem}

\begin{proof}
Let 
$
xb^{-1}=v_0+v_1(cb^2)+\ldots + v_{n-1}(cb^2)^{n-1}
$ 
be the presentation
of the element $xb^{-1}$ in the basis $\Bb(cb^2)$.
Our goal is to show that the subset of the elements 
$b\in V_c$ with $v_{n-1}\neq 0$ 
is non-empty and open in $S^*$.

Multiplying the presentation of $xb^{-1}$ by $b$ we get the following equation:
\begin{equation}\tag{*}
x=v_0b+v_1(cb^3)+\ldots + v_{n-1}c^{n-1}b^{2n-1}.
\end{equation}
Now consider the primitive element $c$ of $S$.
Let $x=x_0+x_1c+\ldots +x_{n-1}c^{n-1}$ 
and 
$b=b_0+b_1c+\ldots+b_{n-1}c^{n-1}$ be the presentations of the
elements $x$ and $b$ in the basis $\Bb(c)$.
Consider the equation (*) in terms of the basis $\Bb(c)$.
We get the system of $n$ linear equations in $n$ variables 
$v_0,\ldots,v_{n-1}$: 
$$
\begin{pmatrix}
x_0\cr
x_1\cr
\vdots\cr
x_{n-1}\cr
\end{pmatrix}
=v_0\cdot 
\begin{pmatrix}
h_{0,0}(b_0,\ldots,b_{n-1})\cr
h_{1,0}(b_0,\ldots,b_{n-1})\cr
\vdots\cr
h_{n-1,0}(b_0,\ldots,b_{n-1})\cr
\end{pmatrix} 
 + \ldots v_{n-1}\cdot 
\begin{pmatrix}
h_{0,n-1}(b_0,\ldots,b_{n-1})\cr
h_{1,n-1}(b_0,\ldots,b_{n-1})\cr
\vdots\cr
h_{n-1,n-1}(b_0,\ldots,b_{n-1})\cr
\end{pmatrix},
$$
where $h_{i,j}$ are polynomials in $n$ variables $b_0,\ldots,b_{n-1}$
with coefficients from $k$.
In particular, we have $h_{i,0}=b_i$ for all $i=0\ldots n-1$.
Solving this linear system we get 
$$
v_{n-1}=\det(A_{(n-1)})/\det(A),
$$
where $A$ is the matrix of the system and the matrix 
$A_{(n-1)}$ is got by replacing
the last column of $A$ by the vector $x$ (of free terms).
Observe that both $\det(A_{(n-1)})$ and $\det(A)$ are homogeneous polynomials
(in $n$ variables $b_0,b_1,\ldots,b_{n-1}$)
of degrees $1+3+\ldots+(2n-3)=(n-1)^2$ and $1+3+\ldots+(2n-1)=n^2$ 
respectively.
Hence, the map $V_c\ra k$ given by $b=(b_0,\ldots,b_{n-1})\mapsto v_{n-1}$
is the regular map (observe that the determinant $\det(A)$ is non-zero
for all $b\in V_c$, since $cb^2$ is primitive). 
Thus, the subset $W_c$ is open in $S^*$.
The fact that $W_c$ is non-empty follows from Sublemma \ref{vanpol} below.
\end{proof}

\begin{sublemma}\label{vanpol}
If $x\neq 0$, then the polynomial $\det(A_{n-1})$ is non-trivial.
\end{sublemma}

\begin{proof}
We have the presentation of the determinant 
$$
\det(A_{(n-1)})= \Delta_{0}x_{n-1}-
\Delta_{1}x_{n-2}+\ldots +(-1)^{n-1}\Delta_{n-1}x_0,
$$
where $\Delta_i$, $i=0,\ldots,n-1$, is the $(n-1-i,n-1)$-minor of $A_{(n-1)}$.
For each monomial $b_0^{m_0}b_1^{m_1}\ldots b_{n-1}^{m_{n-1}}$
we define the weight to be the sum $\sum_{i=0}^{n-1} i\cdot m_i$.
Observe that each minor $\Delta_i$ consists of monomials of 
weight $\geq i$.

Now each minor can be viewed
as the sum $\Delta_i= h_i+ h_{>i}$ of two homogeneous polynomials 
(in $n$-variables $b_0,b_1,\ldots,b_{n-1}$) of degree $(n-1)^2$,
where $h_i$ is the sum of monomials of weight $i$ and $h_{>i}$
is the sum of monomials of weight strictly bigger than $i$.
The point is that the coefficients of the polynomial $h_i$ don't
depend on the coefficients of the minimal polynomial $p_c$ for $c$.

We claim that $h_i\neq 0$.
Indeed, for $i=0$ we have $h_0=b_0^{(n-1)^2}$, i.e.,
$h_0\neq 0$.
For $i=1$ we have $h_1= -C^1_{2n-3} b_0^{n(n-2)}b_1=-(2n-3)b_0^{n(n-2)}b_1$
which is non-zero if the characteristic of $k$ doesn't divide $2n-3$.
For $i>1$ the polynomial $h_i$ contains the unique monomial (the monomial
with the maximal power of $b_0$) 
$(-1)^i(2(n-i)-1)b_0^{n(n-2)}b_i$ which is non-zero if the characteristic
of $k$ doesn't divide $2(n-i)-1$.

Now it follows immediately that the polynomials $\Delta_i$, $i=0,\ldots,n-1$,
are linearly independent. And we are done.
\end{proof}
 
We will need one more fact concerning the polynomials $\Delta_i$:
\begin{sublemma}\label{vanprod}
For any integers $0\leq i,j \leq n-1$
The polynomials $\Delta_i\Delta_j$, where $i\leq j$,
are linearly independent.
\end{sublemma}

\begin{proof}
The product $\Delta_i\Delta_j$ is the homogeneous polynomial of
degree $2(n-1)^2$ and can be represented as follows
$\Delta_i\Delta_j=h_ih_j+g_{>i+j}$, where $h_ih_j$ is the sum
of monomials of weight $i+j$ and $g_{>i+j}$ consists of monomials
of weight strictly bigger than $i+j$.
Observe now that $h_ih_j$ contains the unique monomial (the monomial
with the maximal power of $b_0$) 
$(-1)^{i+j}(2(n-i)-1)(2(n-j)-1)b_0^{2n(n-2)}b_ib_j$ 
(see the proof of the previous Sublemma). 
\end{proof}

Now we are ready to prove the main result of this section:
\begin{thm}\label{genpos}
Let $S$ be a simple extension of degree $n$ of a field $k$ of characteristic 
$0$. Let $c$ be a primitive element of $S/k$.
Let $q$ be a regular quadratic form over $k$ of rank $m$.
and 
$x=(x^{(1)},\ldots,x^{(m)})$ be a vector in $S^m$ such that
$q(x)\neq 0$.
By $\{x,c\}=(x_{n-1}^{(1)},\ldots,x_{n-1}^{(m)})$ 
we denote the vector of the $(n-1)$-th  coordinates
of $x$ in the basis $\Bb(c)$, i.e., $\{x,c\}_j=\{x^{(j)},c\}$.
Then the subset 
$$
U_{c,x,q}=\{b\in V_c | q(\{xb^{-1},cb^2\})\neq 0\}
$$
is non-empty and open in $S^*$.
\end{thm}

\begin{rem} According to \ref{charrem} the Theorem is not true
if the characteristic of the residue field $k$ is non-zero.
Take $k=\bar{{\mathbb F}}_3$, $S=k[t]/(t^3-1)$, $q(x)=x^2$ and $c=x=t$.
\end{rem}

\begin{proof}
The proof is a little modification of the proof of \ref{cxprim}.
We use the notation introduced in the proof of Lemma \ref{cxprim}.

Clearly, $U_{c,x,q}$ is open (by the same arguments as in \ref{cxprim}).
The main problem is to show that $U_{c,x,q}$ is non-empty.
Hence, we have to prove that the polynomial 
$q(\det(A_{(n-1)}^{(1)}),\ldots,\det(A_{(n-1)}^{(m)}))$ is non-trivial,
where $A_{(n-1)}^{(j)}$ denotes the matrix corresponding to the element
$x^{(j)}$. 
Let 
$$\det(A_{(n-1)}^{(j)})= \Delta_{0}x_{n-1}^{(j)}-
\Delta_{1}x_{n-2}^{(j)}+\ldots +(-1)^{n-1}\Delta_{n-1}x_0^{(j)}
$$
be the representation as in the proof of \ref{vanpol}, where
$x^{(j)}=\sum_{i=0}^{n-1} x_i^{(j)}c^i$ is the presentation of $x^{(j)}$
in the basis $\Bb(c)$.
Let $q(x)=\sum_j a_j(x^{(j)})^2$ be our quadratic form.

Then, we have
$$
(\det A)^2 \cdot q(\{xb^{-1},cb^2\})=\sum_j a_j(\det A_{(n-1)}^{(j)})^2=
$$
$$
=\sum_j a_j(\sum_i (-1)^ix_{n-1-i}^{(j)}\Delta_i)^2
$$
Now if we replace $\Delta_i$ by $(-1)^it^{n-1-i}$ we get precisely
\begin{equation}
\tag{*}
\sum_j a_j(\sum_i (-1)^ix_{n-1-i}^{(j)}\Delta_i)^2=
q(x^{(1)}(t),\ldots,x^{(m)}(t))
\end{equation}
as the polynomial in $R[t]$, where $x^{(j)}(t)=\sum_i x_i^{(j)}t^i$.

Assume that the polynomial 
$q(\{xb^{-1},cb^2\})$ is trivial.
Since the polynomials $\Delta_i\Delta_j$ are linearly independent
(by Sublemma \ref{vanprod}), this implies that
the polynomials $\Delta_i\Delta_j$ in the sum (*) have zero coefficients.
In particular, the $t^i$ have trivial coefficients as well,
i.e., the polynomial $q(x^{(1)}(t),\ldots,x^{(m)}(t))$ is trivial.
This contradicts with the hypothesis of the Theorem that the image
of $q(x^{(1)}(t),\ldots,x^{(m)}(t))$ by means of the canonical map
$R[t]\ra R[t]/p_c(t)=S$, i.e., precisely $q(x)$, is non-trivial.
\end{proof}

\section{The Knebusch's Norm Principle}

\begin{ntt}
Let $q$ be a regular quadratic form over $R$ of rank $m$.
We denote by $q_S=q\otimes_R S$ the base change form of $q$.
Thus, $q_S$ is the quadratic form over $S$ whose coefficients
come from the base ring $R$.
Let $D(q_S)$ be the subgroup of $S^*$ generated by all invertible
elements of $S$ that are represented by the quadratic form $q_S$, i.e.,
$D(q_S)=\{b\in S^*| b=q_S(x),\; x\in S^m\}$.
By $D(q)$ we denote the respective subgroup 
of values of the form $q$
over $R$.
\end{ntt}

The goal is to prove the following theorem

\begin{thm}\label{normpr} Let $R$ be a local domain
with the residue field $k$ of characteristic 0.  
Let $S/R$ be a simple extension
of degree $n$. Let $q$ be a regular quadratic form over $R$ of rank $m$. 
Then we have the following inclusion between the subgroups of $S^*$:
$$
N^S_R(D(q_S)) \subset D(q),
$$
where $N^S_R$ is the norm map for the extension $S/R$ and
$D(q_S)$ is the group generated by values of the form $q$ over $S$. 
\end{thm}

\begin{lemma}\label{squares}
In the hypothesis of Theorem \ref{normpr} we have
$(S^*)^2\subset D(q_S)$.
\end{lemma}

\begin{proof}
Let $r\in R^*$ be a value of the quadratic form $q$, i.e., $r=q(y)$
for some $y\in R^m$.
Then for any $b\in S^*$ we have $b^2=q(yb)q(y/r)\in D(q_S)$.
\end{proof}

\begin{proof}[The Proof of Theorem \ref{normpr}]
We prove by induction on the degree $n$ of the simple extension $S/R$.
The case $n=1$ is trivial.
Assume $n>1$.

Let $c=q_S(x)$ for some $x\in S^m$. We want to show $N^S_R(c)\in D(q)$.
By \ref{primelts} and Lemma \ref{cprim}
the element $c$ can be written as the product $c=cb^2\cdot (1/b)^2$, 
where $cb^2=q_S(xb)\in D(q_S)$ is primitive.
Hence, by multiplicativity of the norm map and Lemma \ref{squares}
we may assume $c$ (replaced by $cb^2$) is primitive.

Now we mimic the step 3 of the proof of \cite[VII, 5.1]{La}.
Replace $c$ by it's inverse $c^{-1}$. We get the equation $1=cq_S(x)$,
where $c$ is primitive.
More precisely, we have
\begin{equation}
\tag{*}
1=cq_S(x^{(1)}(c),\ldots,x^{(m)}(c)),
\end{equation}
where $x^{(j)}(c)\in S$ is the $j$-th coordinate of the vector $x$ written
in the basis $\Bb(c)$. 
According to Theorem \ref{genpos}
we may assume that the value of the quadratic form $q_S$
on the last coefficients of the vectors $x^{(j)}$, i.e., 
$q(x^{(1)}_{n-1},\ldots,x^{(m)}_{n-1})$, is invertible in $R$.
In fact, it is enough to consider the quotient modulo the maximal
ideal of $R$ (see \ref{primelts}), i.e., the simple extension $\bar{S}/k$.
The open subset $U_{c,x,q}$ from \ref{genpos} is non-empty
and open. Take any element $b$ from $U_{c,x,q}$
and replace $c$ by $cb^2$ and $x$ by $x/b$.  

Consider the pull-back of the equation (*) by means of the canonical map
$R[t]\ra R[t]/p_c(t)=S$.
Since $tq(x(t))-1$ lies in the principal ideal $(p_c(t))$ of 
the polynomial ring $R[t]$ there is
a polynomial $h(t)$ such that 
\begin{equation}
\tag{**}
1+p(t)h(t)=tq_S(x^{(1)}(t),\ldots,x^{(m)}(t)).
\end{equation}
Since $R$ is a domain
the leading coefficient of the left hand side of (**)
coincides with the leading coefficient of $h$, denoted by $r$, 
and coincides with the leading
coefficient of the right hand side that is
$r=q(x^{(1)}_{n-1},\ldots,x^{(m)}_{n-1})$,
where $r$ is invertible in $R$.
Clearly, we have $n+\deg h=2(n-1)+1$. So that $\deg h=n-1$.

As in \cite{La} we have 
$N^S_R(c)=(-1)^{n-1}p_c(0)=(-1)^n/r\cdot (h(0)/r)^{-1}$, where
$r\in D(q)$ and $h(0)/r$ is the norm (up to sign) 
of the respective primitive element $u=t$
of the simple extension $T=R[t]/g(t)$,
where $g(t)=h(t)/r$ is the monic polynomial of degree $\deg h=n-1$.
Observe that taking (**) modulo the principal ideal (g(t))
we get similar to (*)
$$
1=uq_T(x(u)\, mod \, g(u)),\quad
\text{i.e.,}\quad 
u\in D(q_T).
$$
Hence, we may apply the induction hypothesis and conclude that 
$g(0)=h(0)/r \in D(q)$.
Then we immediately get $N^S_R(c)\in D(q)$.
\end{proof}

\section{Grothendieck-Serre's conjecture for the case of spinor group}

\begin{ntt}
Let $R$ be a local domain with the residue field $k$
of characteristic different form 2 
and $q$ be a regular quadratic form over $R$.
Following \cite{Kn} we define the spinor group (scheme) $Spin_q$
to be $Spin_q(R)=\{x\in S\Gamma_q(R)| x\sigma(x)=1\}$, 
where $\sigma$ is the canonical involution. Recall that
$S\Gamma_q(R)=\{c\in C_0(V,q)^* | cVc^{-1}\subset V\}$, 
where $C_0(V,q)$ is the even part of the Clifford algebra of
the respective quadratic space $(V,q)$ over $R$.
\end{ntt}

The goal of the present section is to show:
\begin{thm}\label{groth}
Let $R$ be a local regular ring of geometric type over the base field
$k$ of characteristic 0. Let $K$ be it's quotient field.
Let $q$ be a regular quadratic form over $R$.
Then the induced map on the sets of principal homogeneous spaces
$$
H^1_{et}(R,Spin_q)\ra H^1_{et}(K,Spin_q)
$$
has trivial kernel, 
where $Spin_q$ is the spinor group for the quadratic form $q$.
\end{thm}

\begin{rem}
Observe that Theorem is the particular case $G=Spin_q$ 
of Grothendieck-Serre's conjecture on principal homogeneous spaces 
\cite{Gr},
which states for a flat reductive group scheme $G$ over $R$
the induced map $H^1_{et}(R,G)\ra H^1_{et}(K,G)$ has trivial kernel.
\end{rem}

\begin{proof}
Let $R$ be a local regular ring and $K$ be it's quotient field.
We have the following commutative digram (see \cite{Kn}):
$$
\xymatrix{
SO_q(R) \ar[r]^{SN} \ar[d] & R^*/(R^*)^2 \ar[r] \ar[d] &
H^1_{et}(R, Spin_q) \ar[r] \ar[d] &  H^1_{et}(R,SO_q) \ar[d] \\
SO_q(K) \ar[r]^{SN} & K^*/(K^*)^2 \ar[r] &
H^1_{et}(K, Spin_q) \ar[r] &  H^1_{et}(K,SO_q),
}
$$
where $SN: SO_q(R) \ra H^1_{et}(R,\mu_2)=R^*/(R^*)^2$ is the spinor norm.
The main result of paper \cite{OP} says that the vertical arrow on the right
hand side has trivial kernel (see also \cite[3.4]{Za}). 
Thus, in order to show that the middle one
has trivial kernel it is enough to check that the induced map on the
cokernels $coker(SN)(R)\ra coker(SN)(K)$ is injective. 

Consider the group scheme
$F:S\mapsto coker(SN)(S)$.
According to \cite[s.2]{Za} to prove the mentioned injectivity
we have to show that the functor $F$ satisfies all the axioms of
\cite[s.1, 2]{Za}. In fact, all the axioms, excluding the existence
of transfer map, holds by the same arguments as in \cite[s.3]{Za}.
So that the only thing we have to check is that
for any finite surjective extension $S/R$ of semi-local rings the norm
map $N^S_R$ commutes with the spinor norm.
In fact, it turns out that we can
consider only finite etale extensions $S/R$ of semi-local rings -- 
just replace
the geometric presentation lemma used in \cite{Za} by it's stronger version
from \cite[Theorem 6.1]{Za1} (the proof will be the same).

Observe that by definition of the spinor norm \cite{Kn} it is the same
as to show that the norm map commutes with the
functor $D: S \mapsto q(S)$ (that sends any $R$-algebra $S$ to the
group of values of the quadratic form $q$ on $S$), i.e.,
$N^S_R(D(S)) \subset D(R)$.
Hence, we have to prove the analog of Knebusch's Norm Principle
for quadratic forms in the case of finite etale extensions
of semi-local rings. 
But this is done in Theorem \ref{normpr}.
\end{proof}

\end{document}